\documentclass[a4paper,reqno,11pt]{amsart}
\usepackage{amsfonts,amsmath,amsthm,amssymb,stmaryrd}
\usepackage{mathrsfs}
\usepackage{cases}
\usepackage{mathtools}
\usepackage{color}
\usepackage[colorlinks,linkcolor=blue,citecolor=blue]{hyperref}

\newtheorem{Theorem}{Theorem}[section]
\newtheorem{Lemma}{Lemma}[section]

\newtheorem{Remark}{Remark}[section]

\numberwithin{equation}{section} \allowdisplaybreaks
\makeatletter

\newcommand{\Rmnum}[1]{\expandafter\@slowromancap\romannumeral #1@}
\makeatother
\usepackage[left=1 in, right=1 in,top=1 in, bottom=1 in]{geometry}
\makeatletter
\@namedef{subjclassname@2020}{\textup{2020} Mathematics Subject Classification}  
\makeatother

\begin{document}
\title[Uniqueness of blowup at singular points for superconductivity problem]{\bf Uniqueness of blowup at singular points for superconductivity problem}
\author{Lili Du}
\address{School of Mathematical Science, Shenzhen University, Shenzhen, 518061, P.R. China. and Department of Mathematics, Sichuan University, Chengdu, 610064, P.R. China.}
\email[L. Du]{dulili@szu.edu.cn}
\author{Xu Tang$^*$}
\address{Department of Mathematics, Sichuan University, Chengdu, 610064, P.R. China.}
\email[X. Tang]{tangxu8988@163.com}
\author{Cong Wang}
\address{Department of Mathematics, Sichuan University, Chengdu, 610064, P.R. China.}
\email[C. Wang]{wc252015@163.com}

\thanks{*Corresponding author: tangxu8988@163.com}
\thanks{
This work is supported by the National Natural Science Foundation of China grant 11971331, 12125102, 12301258, and Sichuan Youth Science and Technology Foundation 2021JDTD0024.}

\begin{abstract}
In this paper, we  prove that the uniqueness of blowup at the maximum point of coincidence set of the superconductivity problem, mainly based on the Weiss-type and Monneau-type monotonicity formulas, and the proof of the main results in this paper is inspired the recent paper \cite{CFL22} by Chen-Feng-Li.
\noindent

{\bf \normalsize Keywords }  {
Superconductivity, Obstacle problem, Singular points, Blowup, Uniqueness

}\bigbreak

\end{abstract}


\maketitle

\section{Introduction}
\setcounter{equation}{0}

In this paper, we consider the  superconductivity problem
\begin{equation}\label{11}
	\Delta u=\chi_{\{|\nabla u|>0\}}  \quad \text{in} ~B_1,
\end{equation}
which is an obstacle-type problem derived from superconducting models, with the more general form
\begin{equation}\label{12}
	\Delta u=f(x, u) \chi_{\{|\nabla u|>0\}} \quad \text { in } B_1,
\end{equation}
where $B_1$ is the unit ball in $\mathbb{R}^n$, the function $f>0$ and $f \in C^{0,1}\left(\mathbb{R}^n \times \mathbb{R}\right)$.

In analyzing the evolution of vortices generated in the mean-field model of a magnetic field passing through a superconductor, we obtain a degenerate parabolic-elliptic system. The simplified stationary model of the problem (in a local setting) is reduced to the problem \eqref{12} with appropriate boundary conditions (see \cite{Cha95}).  Berestycki-Bonnet-Chapman \cite{BBC94} and Chapman-Rubinstein-Schatzman  \cite{CRS96} proposed a related model, with a rigorous derivation from the Ginzburg-Landau model by Sandier-Serfaty \cite{SS00}. We would like to refer the interested readers to the references \cite{PSU12} and \cite{Ro} for the physical background.

From the structure of the equation, the problem \eqref{11} is more general than the no-sign obstacle problem, because the function $u$ may take different constant values in different connected branches of  {\it the coincidence set} $\{|\nabla u| =0\}$, which also leads to the complexity of the free boundary of the problem. Elliott-Schatzle-Stoth \cite{ESS98} studied the above general degenerate parabolic-elliptic system, they proved the existence and uniqueness of the viscosity solution in two dimensions, and found the special solutions of the stationary problem. Caffarelli-Salazar \cite{CS02} constructed the viscosity solution of a fully nonlinear elliptic equation more general than problem \eqref{12} and obtained some properties of the viscosity solution. Besides, based on the results in \cite{CS02}, the existence and regularity can be proved using the Alexandroff-Bakelman approximation technique  appropriately. Bonnet-Monneau \cite{BM00} and Monneau \cite{Mon04} investigated the free boundary of a specific configuration (with single patches). To be specific, Bonnet-Monneau \cite{BM00} shows that the existence and regularity via Nash-Morse theory, and  Monneau \cite{Mon04} proved the regularity of the free boundary when it is close enough to the fixed boundary, they also gave a result of stability of the free boundary and gave a bound on the Hausdorff measure of the free boundary. The free boundary in general was first studied by Caffarelli-Salazar \cite{CS02} and then by Caffarelli-Salazar-Shahgholian \cite{CSS04}. In particular, based on  a refined analysis,  Caffarelli-Salazar-Shahgholian  reduced the problem to the one-patch case, and whose global solutions  was  characterized  with the help of  Weiss's monotonicity formula. In \cite{PSU12}, the authors systematically studied the problem \eqref{11}, and proved the optimal regularity and nondegeneracy of the solution. Moreover, they established the free boundary regularity near the regular point and the structure of the regular point set by blowup method. However, as mentioned in \cite[Chapter 7, Notes]{PSU12}, for the problem \eqref{11}, very little is known about the singular set, all the existing methods seem to fail.
This paper attempts to make a preliminary analysis of the singular set of the problem \eqref{11}. It is well know the blowup limit  possible dependence on subsequences is one of the main diffculties in the study of free boundary problems. This paper is devoted to the uniqueness of blowup at singular points for the superconductivity problem \eqref{11}.

For the classical obstacle problem, the structure of the singular set  was discovered by Caffarelli \cite{Caf98}. Subsequently, Monneau \cite{Mon03} introduced a concise method to prove the uniqueness of blowup at singular points, and the essential tool is the Monneau's monotonicity formula, on which one can prove that the continuous dependence of blowups, and study the structure of the singular set (see \cite[Chapter 7.4]{PSU12}). Recently, Chen-Feng-Li \cite{CFL22} noted that by means of the Monneau's monotonicity formula it is also possible to prove that the uniqueness of blowup at singular points of the no-sign obstacle problem.  Moreover, this paper attempts to generalize the Monneau's monotonicity formula to the superconductivity problem \eqref{11}.

Throughout this paper, we denote by
$\Gamma:=\partial\{ |\nabla u|>0\} \cap B_1$ the free boundary for the problem \eqref{11}, and $\Sigma$ {\it the singular set} of $\Gamma$, i.e. $x^{0} \in \Sigma$ if and only if there exists a sequence $r_j \rightarrow 0$ such that $u_{x_0, r_j}(x):=\frac{u\left(x^0+r_j x\right)-u(x^0)}{r_j^2}$ converges to a homogeneous quadratic polynomial $q(x)$ with $\Delta q=1$. We call such $q(x)$ is  a {\it blowup} of $u$ at $x^0$. Without loss of generality, let us set $x^0=0$ below. For convenience, we denote
\begin{align*}
	\mathcal{Q}:=&\{q(x) \text { homogeneous quadratic polynomial}: \Delta q=1\}
\end{align*}
and
\begin{align*}
	\mathcal{Q}^{+}:=&\{q \in \mathcal{Q}: q \geq 0\}.
\end{align*}

Our main result reads

\begin{Theorem}\label{Thm11}
Let $u$ be any solution to the problem \eqref{11}, and that $0 \in \Sigma$ satisfies
\begin{align}\label{13}
	u \leq u(0) \quad \text{on} \quad \{|\nabla u|=0\}.
\end{align}
Then there is a  $q_0(x) \in \mathcal{Q}$ such that
\begin{align*}
	u_r(x):=\frac{u(r x)-u(0)}{r^2} \rightarrow q_0 ~~\text { in } C_{\rm{loc}}^{1,\alpha}\left(\mathbb{R}^n\right) \quad \text { as } r \rightarrow 0+
\end{align*}
for any $\alpha \in (0,1)$. Moreover, there holds
\begin{align}\label{15}
	u(x)-u(0)=q_{0}(x)+o\left(|x|^2\right).
\end{align}
\end{Theorem}

\begin{Remark}
It should be noted that $u$ be a solution to the problem \eqref{11}  means that $u \in C_{\rm{loc}}^{1,1}\left(B_1\right)$ satisfying $\left\|D^2 u\right\|_{ L^{\infty}\left(B_1\right)} \leq M$ for some constant $M>0$ and $0 \in \Gamma$ (see \cite[Definition 3.15]{PSU12}).
\end{Remark}

\begin{Remark}
 Whether blowup is unique at a singular point is usually a priori unknown, i.e. maybe $u_r$ will sub-converge to a different polynominal for another sequence $\tilde{r}_j \rightarrow 0$. In Theorem \ref{Thm11}, there is no need to take a subsequence for $r \rightarrow 0$, indicating the uniqueness of blowup.
\end{Remark}

Our idea to prove the uniqueness of blowup at  singular point is to construct a corresponding Monneau-type monotonicity formula for the problem \eqref{11}. It is worth noting that a common point between the classical obstacle problem and the no-sign obstacle problem is that both the solution $u$ and its first partial derivative are equal to $0$ at the free boundary point. However, in the problem \eqref{11}, only the gradient of $u$ is zero on the free boundary, and there is no information about the value of $u$ itself, so the Monneau's monotonicity formula of the classical obstacle problem (e.g. see \cite[Theorem 7.4]{PSU12}) is not valid for the problem \eqref{11}.

Inspired by the research on the superconductivity problem in \cite{PSU12}, when we derive the monotonicity formula in next section, we add the hypothesis \eqref{13}.  Furthermore, we consider
\begin{align*}
	\frac{1}{r^{n+3}} \int_{\partial B r}(u-u(0)-q)^2 d \mathcal{H}^{n-1}
\end{align*}
as the Monneau's energy functional, where $q \in \mathcal{Q}^{+}$. We will show that the derivative  of the functional with respect to $r$ (see Lemma \ref{Lem23} below for details)
\begin{align*}
	\frac{2}{r^{n+4}} \int_{\partial B r} w(x \cdot \nabla w-2 w) d \mathcal{H}^{n-1}
\end{align*}
is nonnegative, where   $w=u-u(0)-q$.  According to the construction of the Weiss's energy functional (see \eqref{21} below), there is the identity
\begin{align*}
	\frac{1}{r^{n+3}} \int_{\partial B_r} w(x \cdot \nabla w-2 w) d \mathcal{H}^{n-1}=W(r, u)-W(0+, u)+\frac{1}{r^{n+2}} \int_{B_r} w \Delta w d x,
\end{align*}
so the problem is transformed into proving that the right hand side of the above equation is nonnegative. The addition of the first two terms is nonnegative  can be obtained by applying the Weiss's monotonicity formula, therefore a key fact is
\begin{align}\label{14}
	w \Delta w \geq 0 \quad \text { in } \quad B_r.
\end{align}
This fact also plays an important role in the discussion of the classical obstacle problem. In addition, it should be noted that the limit $W(0+, u)$ is required to exist in our proof, so monotone nondecreasing of the Weiss's energy functional is a sufficient condition. The hypothesis \eqref{13} guarantees that the Weiss's monotonicity formula  and \eqref{14} holds simultaneously.

Since the blowup of the problem \eqref{11} is known to be a homogeneous quadratic polynomial in $\mathcal{Q}$ at every singular point, in order to prove the uniqueness of $q_0 \in \mathcal{Q}$ in Theorem \ref{Thm11}, it is only necessary to prove that the corresponding coefficient matrix is unique. It is worth noting that $q \in \mathcal{Q}^{+}$ is arbitrary. We use the method similar to that in \cite{CFL22}, select a special family of $\{q^{t}\}_{t \in (-1,1)}$, where $q^{t} \in \mathcal{Q}^{+}$. Applying the Monneau's monotonicity formula, we get an equation with parameter $t$. The derivative of both sides of this equation with respect to $t$ still yields an identity. Subtly, the latter leads to the desired uniqueness.

In addition, \eqref{15} such estimates are the starting point for studying the structure of the singular set in problem \eqref{11}. If we can show that the bound $o(r^2)$ on the right hand side of \eqref{15} is uniform for all singular point, or derive a more accurate quantitative characterization of $o(r^2)$, we can expect to achieve the corresponding results for the structure of the singular set to be somewhat similar to those for the classical obstacle problem (see \cite{FR22, Fig18, Fig18b, PSU12} and the reference therein).

The structure of this paper is arranged as follows. In Section \ref{secion2}, we derive the monotonicity formulas.  In Section \ref{section3}, we prove the main result of this paper.

\section{Monotonicity formulas}\label{secion2}

In the present section, we first introduce the Weiss's monotonicity formula and then derive the Monneau's monotonicity formula from it.

\begin{Lemma}[Weiss's monotonicity formula]\label{Lem21}
Let $u$ be any solution to the problem \eqref{11}, and that
\begin{align}\label{25}
u \leq u(0) \quad \text{on} \quad \{|\nabla u|=0\}.
\end{align}
Then
\begin{align}\label{21}
	r \mapsto W(r, u):=\frac{1}{r^{n+2}} \int_{B_r}\left(|\nabla u|^2+2(u-u(0))\right) d x-\frac{2}{r^{n+3}} \int_{\partial B_r}(u-u(0))^2 d\mathcal{H}^{n-1}
\end{align}
is a nondecreasing absolutely continuous function for $0<r<1$ and
\begin{align}\label{22}
	\frac{d}{d r} W\left(r, u\right) \geq \frac{2}{r^{n+4}} \int_{\partial B_r}\left|x \cdot \nabla u-2(u-u(0))\right|^2 d \mathcal{H}^{n-1}
\end{align}
for a.e. $0<r<1$.
\end{Lemma}

\begin{Remark}
It is worth mentioning that the functional $W$ has the following scaling property
\begin{align}\label{23}
	W\left(r s, u\right)=W(s, u_r)
\end{align}
for any $0<r<1,~0<s<\frac{1}{r}$, where
\begin{align}\label{24}
	u_r(x)=\frac{u(r x)-u(0)}{r^2} .
\end{align}
In particular,
\begin{align}\label{212}
	W\left(r, u\right)=W(1, u_r).
\end{align}
\end{Remark}

The proof can refer to  \cite[Theorem 3.26]{PSU12}, but for completeness, the proof is also given below.

\noindent{\it{Proof of Lemma \ref{Lem21}.}}
It follows from \eqref{212}, we have
\begin{align*}
\frac{d}{d r} W(r, u)&=\frac{d}{d r} W\left(1, u_r\right)\\
		& =\int_{B_1}\frac{d}{d r}\left(\left|\nabla u_r\right|^2+2 u_r\right) d x-2 \int_{\partial B_1} \frac{d}{d r}\left(u_r^2\right) d \mathcal{H}^{n-1}  \\
		& =\int_{B_1}\left(2 \nabla u_r \cdot \nabla \frac{d u_r}{d r}+2 \frac{d u_r}{d r}\right) d x-4 \int_{\partial B_1} u_r \frac{d u_r}{d r} d\mathcal{H}^{n-1} .
\end{align*}
By integrating by parts, we get
\begin{align*}
	\int_{B_1} \nabla u_r \cdot \nabla \frac{d u_r}{d r} d x= \int_{\partial B_1} \frac{d u_r}{d r} \frac{\partial u_r}{\partial \nu} d\mathcal{H}^{n-1}-\int_{B_1} \Delta u_r \frac{d u_r}{d r} d x,
\end{align*}
where  $\frac{\partial u_r}{\partial \nu}$ is the outer normal derivative of $u_r$ on $\partial B_1$, and so
\begin{align}\label{219}
		\frac{d}{d r} W(r, u)
		=2 \int_{B_1} \frac{d u_r}{d r}\left(1-\Delta u_r\right) d x+2 \int_{\partial B_1} \frac{d u_r}{d r}\left(\frac{\partial u_r}{\partial \nu}-2 u_r\right) d \mathcal{H}^{n-1}.
\end{align}

Note that
\begin{align*}
1-\Delta u_r=\chi_{\left\{\left|\nabla u_r\right|=0\right\}}
\end{align*}
and
\begin{align}
		\frac{d u_r}{d r}
		& =\frac{r x \cdot \nabla u(r x)-2(u(rx)-u(0))}{r^3} \notag\\
		& =\frac{x}{r} \cdot \nabla u_r-\frac{2}{r} u_r \label{26}
\end{align}
implies that
\begin{align*}
	\frac{d u_r}{d r}=-\frac{2}{r} u_r \quad \text { on }\quad \left\{\left|\nabla u_r\right|=0\right\} .
\end{align*}
Hence, we obtain
\begin{align}
		\int_{B_1} \frac{d u_r}{d r}\left(1-\Delta u_r\right) d x & =-\frac{2}{r} \int_{\{|\nabla u_r|=0\}} u_r d x \notag\\
		& =-\frac{2}{r} \int_{\left\{\left|\nabla u_r\right|=0\right\}} \frac{u(r x)-u(0)}{r^2} d x \notag\\
		& \geq 0, \label{220}
\end{align}
since the assumption \eqref{25} leads to $u(rx) \leq u(0)$ on $\left\{\left|\nabla u_r\right|=0\right\}$. It follows from \eqref{219} and \eqref{220} that
\begin{align*}
		\frac{d}{d r} W(r, u) & \geq 2 \int_{\partial B_1} \frac{d u_r}{d r}\left(\frac{\partial u_r}{\partial \nu}-2 u_r\right) d \mathcal{H}^{n-1} \\
		& =2 \int_{\partial B_1} \frac{d u_r}{d r}\left(x \cdot \nabla u_r-2 u_r\right) d \mathcal{H}^{n-1} \\
		& =\frac{2}{r} \int_{\partial B_1}\left(x \cdot \nabla u_r-2 u_r\right)^2 d \mathcal{H}^{n-1}  \\
		& =\frac{2}{r^{n+4}} \int_{\partial B_r} (x \cdot \nabla u-2(u-u(0)))^2 d \mathcal{H}^{n-1},
\end{align*}
here we have used the fact \eqref{26}.
\qed

The following results are needed to derive the  Monneau's monotonicity formula.

\begin{Lemma}\label{Lem22}
Assume that $u_{r_j} \rightarrow u_0$ in $C_{\rm{loc}}^{1,\alpha}\left(\mathbb{R}^n\right)$ for some sequence $r_j \rightarrow 0$, then
\begin{align}\label{27}
	W\left(r, u_0\right)=W(r, q)
\end{align}
for any $q \in \mathcal{Q}$ and any $r \in(0,1)$.
\end{Lemma}

\begin{proof}
Noting that $0$ is a singular point for the problem \eqref{11}, we know that $u_0$ is a $2$-homogeneous polynomial (see \cite[Theorem 3.23]{PSU12}), i.e. $u_0 \in \mathcal{Q}$.

Thanks to Lemma \ref{Lem21} and the scaling property \eqref{23}, we obtain
\begin{align}\label{28}
	W\left(r, u_0\right)=\lim _{j \rightarrow \infty} W\left(r, u_{r_j}\right)=\lim _{j \rightarrow \infty} W\left(r r_j, u\right)=W(0+, u)
\end{align}
for any $r>0$, which gives that $W\left(r, u_0\right)$ is constant. In particular,
\begin{align}\label{29}
		W\left(r, u_0\right)  \equiv W\left(1, u_0\right)
		 =W(0+, u) .
\end{align}
Taking $r=1$ in \eqref{28}, we get
\begin{align*}
		W(0+, u) & =W\left(1, u_0\right) \\
		& =\int_{B_1}\left(\left|\nabla u_0\right|^2+2\left(u_0-u_0(0)\right)\right) d x-2 \int_{\partial B_1}\left(u_0-u_0(0)\right)^2 d \mathcal{H}^{n-1} \\
		& =\int_{B_1}\left(\left|\nabla u_0\right|^2+2 u_0\right) d x-2 \int_{\partial B_1} u_0^2  d \mathcal{H}^{n-1} \\
		& =\int_{\partial B_1} u_0 \frac{\partial u_0}{\partial \nu} d \mathcal{H}^{n-1}-\int_{B_1} u_0 \Delta u_0 d x+2 \int_{B_1} u_0 d x-2 \int_{\partial B_1} u_0^2 d \mathcal{H}^{n-1} \\
		& =\int_{B_1}\left(-\Delta u_0+2\right) u_0 d x+\int_{\partial B_1} u_0\left(\frac{\partial u_0}{\partial \nu}-2 u_0\right) d \mathcal{H}^{n-1}.
\end{align*}
In virtue of $u_0 \in \mathcal{Q}$, then
 $$\frac{\partial u_0}{\partial \nu}-2 u_0=x \cdot \nabla u_0-2 u_0=0 \quad \text{on}\quad \partial B_1.$$
In addition,
\begin{align*}
	\Delta u_0=\chi_{\left\{\left|\nabla u_0\right|>0\right\}}.
\end{align*}
Hence,
\begin{align*}
		W(0+, u) & =\int_{B_1}\left(-\Delta u_0+2\right) u_0 d x \\
		& =\int_{B_1 \cap\left\{\left|\nabla u_0\right|>0\right\}}\left(-\Delta u_0+2\right) u_0 d x+\int_{B_1 \cap\left\{\left|\nabla u_0\right|=0\right\}}\left(-\Delta u_0+2\right) u_0 d x \\
		& =\int_{B_1 \cap\left\{\left|\nabla u_0\right|>0\right\}} u_0 d x+\int_{B_1 \cap\left\{\left|\nabla u_0\right|=0\right\}} 2 u_0 d x \\
		& =\int_{B_1} u_0 d x+\int_{B_1 \cap\left\{\left|\nabla u_0\right|=0\right\}} u_0 d x \\
		& =\int_{B_1} u_0 d x,
\end{align*}
where the last equality follows from the fact that the set $B_1 \cap\{|\nabla u_0|=0\}$ is measure zero.

Now let $p=\frac{1}{2} x \cdot A x \in \mathcal{Q}$, then
\begin{align}\label{211}
 W(r, p)=W(1, p)
\end{align}
for any $r \in (0,1)$, according to the scaling property.

Next we compute
\begin{align*}
		W(1, p) & =\int_{B_1}\left(|\nabla p|^2+2 p\right) d x-2 \int_{\partial B_1} p^2 d \mathcal{H}^{n-1} \\
		& =\left(\int_{\partial B_1} p \frac{\partial p}{\partial \nu} d \mathcal{H}^{n-1}-\int_{B_1} p \Delta p d x\right)+2 \int_{B_1} p d x-2 \int_{\partial B_1} p^2 d \mathcal{H}^{n-1} \\
		& =\int_{B_1} p d x+\int_{\partial B_1} p (x \cdot \nabla p-2 p) d \mathcal{H}^{n-1} \\
		& =\int_{B_1} p d x .
\end{align*}

Finally, a direct computation shows that there exists a dimensional constant $\alpha_n>0$ such that
\begin{align*}
	\int_{B_1} p d x=\int_{B_1} u_0 d x=\alpha_n.
\end{align*}
Thus, we conclude that \eqref{27}.
\end{proof}

\begin{Lemma}[Monneau's monotonicity formula]\label{Lem23}
Let $u$ be any solution to the problem \eqref{11}, and that $0 \in \Sigma$ satisfies $u \leq u(0)$ on $\{|\nabla u|=0\}$. Then for any $q \in \mathcal{Q}^{+}$, the functional
\begin{align}\label{213}
	r \mapsto M(r, u, q):=\frac{1}{r^{n+3}} \int_{\partial B r}(u-u(0)-q)^2 d \mathcal{H}^{n-1}
\end{align}
is monotone nondecreasing for $r \in(0,1)$.
\end{Lemma}

\begin{Remark}
The functional $M$ has the following nice rescaling property
\begin{align}\label{214}
	M(r, u, q)=M\left(1, u_r, q\right) \quad \text{for all} \quad q \in \mathcal{Q}^{+}.
\end{align}
\end{Remark}

\begin{Remark}
The polynomial $q$ in lemma \ref{Lem23} may not be any blowup limit for the problem \eqref{11}.
\end{Remark}

\noindent\emph{Proof of Lemma \ref{Lem23}.}
Let $w=u-u(0)-q$, then
\begin{align}
		\frac{d}{d r} M(r, u, q) & =\frac{d}{d r}\left(\frac{1}{r^{n+3}} \int_{\partial B r} w^2(x) d \mathcal{H}^{n-1}\right) \notag\\
		& =\frac{d}{d r} \int_{\partial B_1} \frac{w^2(r y)}{r^4} d \mathcal{H}^{n-1} \notag\\
		& =\int_{\partial B_1} \frac{2 w(r y)(r y \cdot \nabla w(r y)-2 w(r y))}{r^5} d \mathcal{H}^{n-1} \notag\\
		& =\frac{2}{r^{n+4}} \int_{\partial B r} w(x \cdot \nabla w-2 w) d \mathcal{H}^{n-1} . \label{215}
\end{align}
By Lemma \ref{Lem22}, we have
\begin{align*}
	W(0+, u)=W(r, q)=\alpha_n .
\end{align*}
Then
\begin{align}
		 W(r, u)-W(0+, u)
		= & W(r, u)-W(r, q) \notag\\
		= & \frac{1}{r^{n+2}} \int_{B r}\left(|\nabla u|^2-|\nabla q|^2+2(u-u(0)-q)\right) d x \notag\\
		& -\frac{2}{r^{n+3}} \int_{\partial B_r}\left((u-u(0))^2-q^2\right) d \mathcal{H}^{n-1} \notag\\
		= & \frac{1}{r^{n+2}} \int_{B r}\left(|\nabla w|^2+2 \nabla w \cdot \nabla q+2 w\right) d x \notag\\
		& -\frac{2}{r^{n+3}} \int_{\partial B r} w(w+2 q) d \mathcal{H}^{n-1}\notag\\
        =&\frac{1}{r^{n+2}} \int_{B r}|\nabla w|^2 d x-\frac{2}{r^{n+3}} \int_{\partial B_r} w^2 d \mathcal{H}^{n-1}\notag\\
        &+\frac{2}{r^{n+3}} \int_{\partial B r} w(x \cdot \nabla q-2 q) d \mathcal{H}^{n-1} \notag\\
		=&\frac{1}{r^{n+2}} \int_{B_r}|\nabla w|^2 d x-\frac{2}{r^{n+3}} \int_{\partial B r} w^2 d \mathcal{H}^{n-1} \notag\\
		=&\frac{1}{r^{n+2}} \int_{B_r}(-w \Delta w) d x+\frac{1}{r^{n+3}} \int_{\partial B_r} w(x \cdot \nabla w-2 w) d \mathcal{H}^{n-1} . \label{216}
\end{align}

On the other hand, we have
\begin{align*}
		w \Delta w
		& =(u-u(0)-q)(\Delta u-1) \\
		&=\left\{\begin{array}{lr}
			0 & \text { on }\{|\nabla u|>0\}, \\
			q-(u-u(0)) & \text { on }\{|\nabla u|=0\} .
		\end{array}\right.
\end{align*}
Since $\Delta u=0$ and $u \leq u(0)$ on $\{|\nabla u|=0\}$, $q \in \mathcal{Q}^{+}$ means that $q \geq 0$, then
\begin{align}\label{217}
	w \Delta w \geq 0 \quad \text { in } ~B_1.
\end{align}
Combining \eqref{216} and \eqref{217}, we arrive at
\begin{align}\label{218}
	W(r, u)-W(0+, u) \leq \frac{1}{r^{n+3}} \int_{\partial B_r} w (x \cdot \nabla w-2 w) d \mathcal{H}^{n-1}.
\end{align}
Thus, it follows from the monotonicity of Weiss's energy functional, and \eqref{215}, \eqref{218} that
\begin{align*}
		\frac{d}{d r} M(r, u, q)  \geq \frac{2}{r}(W(r, u)-W(0+, u))
		 \geq 0 .
\end{align*}
\qed

\section{Proof of the main result}\label{section3}

With the previous preparation, we can now prove the main result of this paper.

\noindent\emph{Proof of Theorem \ref{Thm11}.}  Let us first prove the uniqueness of blowup. Assume that
\begin{align}\label{31}
	u_{r_j} \rightarrow q, ~u_{\tilde{r}_j} \rightarrow \tilde{q} \quad \text { in } C_{\rm{loc}}^{1, \alpha}\left(\mathbb{R}^n\right),
\end{align}
for two sequences $r_j \rightarrow 0, \tilde{r}_j \rightarrow 0$, where $q=\frac{1}{2} x \cdot A x \in \mathcal{Q}$ and $\tilde{q}=\frac{1}{2} x \cdot \tilde{A} x \in \mathcal{Q}$. Note that $A$ and $\tilde{A}$ are two symmetric matrices with $\operatorname{Tr}(A)=\operatorname{Tr}(\tilde{A})=1$. We just need to prove that $A=\tilde{A}$.

Given any $\bar{q}=\frac{1}{2} x \cdot B x \in \mathcal{Q}^{+}$, Lemma \ref{Lem23} gives that
\begin{align}\label{32}
	M\left(1, u_r, \bar{q}\right)=\int_{\partial B_1}\left(u_r-\bar{q}\right)^2 d \mathcal{H}^{n-1} \text { is monotone nondecreasing for } r \in(0,1) .
\end{align}
Moreover, we have that
\begin{align}\label{33}
	\int_{\partial B_1}\left(\frac{1}{2} x \cdot(A-B) x\right)^2 d \mathcal{H}^{n-1}=\int_{\partial B_1}\left(\frac{1}{2} x \cdot(\tilde{A}-B) x\right)^2 d \mathcal{H}^{n-1} .
\end{align}
Since $A-\tilde{A}$ is symmetric, by rotating the coordinates, we may assume that $A-\tilde{A}$ is diagonalized with the eigenvalues  $\lambda_1, \cdots, \lambda_n$. Now, we choose $B=B^t=\left(b_{i j}^t\right)_{i, j=1}^n$, where $b_{11}^t=\frac{1}{2}(1-t), b_{22}^t=$ $\frac{1}{2}(1+t)$, $t \in(-1,1)$, and $b_{i j}^t=0$ for all other $i, j$. It is easy to verify that $B^t \geq 0$ and $\operatorname{Tr}\left(B^t\right)=1$, and then  $\frac{1}{2} x \cdot B^t x \in \mathcal{Q}^{+}$.

Now, recalling the identity \eqref{33} gives that
\begin{align}\label{34}
	f(t):=\int_{\partial B_1}(x \cdot(A-\tilde{A}) x)\left(x \cdot\left(A+\tilde{A}-2 B^t\right) x\right) d \mathcal{H}^{n-1}=0 .
\end{align}
Note that
\begin{align}\label{35}
	x \cdot(A-\tilde{A}) x=\sum_{i=1}^n \lambda_i x_i^2
\end{align}
and
\begin{align}\label{38}
	\operatorname{Tr}(A-\tilde{A})=\sum_{i=1}^n \lambda_i=0.
\end{align}
Let $A=(a_{ij}), \tilde{A}=(\tilde{a}_{ij})$, it  can be calculated directly that
\begin{align*}
	x \cdot\left(A+\widetilde{A}-2 B^t\right) x=\left(a_{11}+\widetilde{a}_{11}-(1+t)\right) x_1^2+\left(a_{22}+\widetilde{a}_{22}-(1-t)\right)x_2^2+R,
\end{align*}
where $R$ is independent on $t$, and then
\begin{align}\label{36}
	\frac{d}{dt}\left(x \cdot\left(A+\tilde{A}-2 B^t\right) x\right)=x_1^2-x_2^2 .
\end{align}
It follows from \eqref{34}, \eqref{35} and \eqref{36} that
\begin{align*}
		0=\frac{d f}{d t} & =\int_{\partial B_1}\left(\sum_{i=1}^n \lambda_i x_i^2\right)\left(x_1^2-x_2^2\right) d \mathcal{H}^{n-1} \\
			& =\int_{\partial B_1}\left(\lambda_1 x_1^2+\lambda_2 x_2^2\right)\left(x_1^2-x_2^2\right) d \mathcal{H}^{n-1} \\
			& \quad+\int_{\partial_{B_1}}\left(\sum_{i=3}^n \lambda_i x_i^2\right)\left(x_1^2-x_2{ }^2\right) d \mathcal{H}^{n-1}\\
		& =\int_{\partial B_1}\left(\lambda_1 x_1^4-\lambda_2 x_2^4-\left(\lambda_1-\lambda_2\right) x_1^2 x_2^2\right) d \mathcal{H}^{n-1} \\
		& \quad +\int_{\partial B_1}\left(\sum_{i=3}^n \lambda_i x_i^2\right)\left(x_1^2-x_2^2\right) d \mathcal{H}^{n-1}.
\end{align*}
Due to the symmetry of $\partial B_1$ we have that
$$\int_{\partial B_1} x_1^4 d \mathcal{H}^{n-1}=\int_{\partial B_1} x_2^4 d \mathcal{H}^{n-1},$$
and that
\begin{align*}
	\int_{\partial B_1} x_1^2 x_i^2 d \mathcal{H}^{n-1}=\int_{\partial B_1} x_2^2 x_i^2 d \mathcal{H}^{n-1} \quad \text { for all } ~~i=3, \cdots, n .
\end{align*}
Hence from the above calculation we have
\begin{align}\label{37}
\frac{d f}{d t}=\left(\lambda_1-\lambda_2\right)\left(\int_{\partial B_1} x_1^4 d \mathcal{H}^{n-1}-\int_{\partial B_1} x_1^2 x_2^2 d \mathcal{H}^{n-1}\right)=0.
\end{align}
We observe that
\begin{align*}
		\int_{\partial B_1} x_1^4 d \mathcal{H}^{n-1} =\frac{1}{4} \int_{\partial B_1} \frac{\partial\left(x_1^4\right)}{\partial \nu} d \mathcal{H}^{n-1}=\frac{1}{4} \int_{ B_1} \Delta\left(x_1^4\right) d x
		 =3 \int_{B_1} x_1^2 d x
\end{align*}
and that
\begin{align*}
		\int_{\partial B_1} x_1^2 x_2^2 d \mathcal{H}^{n-1} & =\frac{1}{4} \int_{\partial B_1} \frac{\partial\left(x_1^2 x_2^2\right)}{\partial \nu} d \mathcal{H}^{n-1}
		=\frac{1}{4} \int_{B_1} \Delta\left(x_1^2 x_2^2\right) d x \\
		& =\frac{1}{2} \int_{B_1}\left(x_1^2+x_2^2\right) d x
		 =\int_{B_1} x_1^2 d x.
\end{align*}
Hence,
\begin{align*}
	\int_{\partial B_1} x_1^4 d \mathcal{H}^{n-1}-\int_{\partial B_1} x_1^2 x_2^2 d \mathcal{H}^{n-1}=2 \int_{B_1} x_1^2 d x>0,
\end{align*}
which implies that $\lambda_1=\lambda_2$ from \eqref{37}.

Similarly, for any $1 \leq i_0<j_0 \leq n$, we can choose $B^t=$ $\left(b_{i j}^t\right)_{i, j=1}^n$, where $b_{i_0 i_0}^t=\frac{1}{2}(1-t), b_{j_0 j_0}^t=\frac{1}{2}(1+t)$, $t \in(-1,1)$, and $b_{i j}^t=0$ for all other $i, j$. From an argument similar to the one above we get  $\lambda_{i_0}=\lambda_{j_0}$. Hence $\lambda_1=\lambda_2=\cdots=\lambda_n$. In virtue of the fact \eqref{38}, we conclude that $\lambda_i=0$ for all $i=1, \cdots, n$. Therefore $A=\tilde{A}$.

Finally, we show that
\begin{align*}
	u(x)-u(0)=q_0(x)+o\left(|x|^2\right) .
\end{align*}
This is equivalent to
\begin{align*}
	r^{-2}\left\|u-u(0)-q_0\right\|_{L^{\infty}\left(B_r\right)} \rightarrow 0 \quad \text { as } r \rightarrow 0 .
\end{align*}
Indeed, assume by contradiction that there is a subsequence $r_k \rightarrow 0$ along which
\begin{align*}
	r_k^{-2}\|u-u(0)-q_0\|_{L^{\infty}\left(B_{r_k}\right)} \geq c
\end{align*}
for some constant $c>0$. Then, there is a subsequence $r_{k_j}$ such that $u_{r_{k_j}} \rightarrow u_0$ in $C_{\mathrm{loc}}^{1,\alpha}\left(\mathbb{R}^n\right)$, for a certain blow-up $u_0$ satisfying
\begin{align*}
	\left\|u_0-q_0\right\|_{L^{\infty}\left(B_1\right)} \geq c.
\end{align*}
On the other hand, the uniqueness of blowup implies that $u_0=q_0$, and hence we reach a contradiction.

 So far, we have completed the proof of Theorem \ref{Thm11}. \qed

\end{document}